\documentclass[12pt]{article}%_if_under_LATeX2e

\usepackage{amsmath, amsthm, amscd, amsfonts, amssymb, graphicx, color}
\usepackage[bookmarksnumbered, plainpages]{hyperref}
\usepackage{amsmath}
\usepackage{amsfonts}
\usepackage{amssymb}

\newtheorem{thm}{Theorem}
\numberwithin{thm}{section}
\newtheorem{lemma} [thm]{{Lemma}}
\newtheorem{proposition} [thm] {Proposition}
\newtheorem{Fact} [thm] {Fact}
\newtheorem{remark}[thm]{Remark}
\newcommand{\Hom}{\mbox{\rm Hom}}%ENDOMORFO
\newcommand{\Aut}{\mbox{\rm Aut}}%AUTOMORFO
\newcommand{\PAut}{\mbox{\rm PAut}}%AUTOMORFISMI_POTENZA
\newcommand{\IAut}{\mbox{\rm IAut}}%
\newcommand{\FAut}{\mbox{\rm FAut}}%AUTOMORFISMI_FINITARI
\newcommand{\St}{\mbox{\rm St}}
\newcommand{\N}{{\Bbb N}}%numeri_naturali_IN_GRASSETTO
\newcommand{\Z}{{\Bbb Z}}%numeri_INTERI
\newcommand{\Q}{{\Bbb Q}}%numeri_rAZIONALI__IN_GRASSETTO
\newcommand{\la}{{\langle}}%parentesi<
\newcommand{\ra}{{\rangle}}%partentesi>
\newcommand{\G}{{\Gamma}}%
\newcommand{\D}{{\Delta}}%

\newcommand{\iso}{\simeq}%accent
\newcommand{\p}{{\varphi}}%VARPHI
\newcommand{\g}{{\gamma}}%
\newcommand{\plus}{\oplus}
\newcommand{\n}{\lhd}%SOTTOGRUPPO_NORMALE
\newcommand{\U}{{\cal U}}%invertibili
\newcommand{\QED}{\hfill $\square$\bigskip}%FINE_DIMOSTRAZIONE_LaTeX2e
\newcommand{\pf}{{\noindent\bf Proof.\ \ }}%PROOF
\title{The group of inertial automorphisms\\ of an abelian group}
\date{\sl \small dedicated  to Martin L. Newell}
\author{ Ulderico Dardano\  -\
Silvana Rinauro }

\begin{document}

\maketitle

\abstract{We study the group $\IAut(A)$ generated by the inertial automorphisms of an abelian group $A$, that is, automorphisms $\g$  with the property that each subgroup $H$ of $A$ has finite index in the subgroup 
generated by $H$ and $H\g$. Clearly, $\IAut(A)$ contains the group $\FAut(A)$ of finitary automorphisms of $A$, which is known to be locally finite.  In a previous paper, we showed that $\IAut(A)$ is (locally finite)-by-abelian. In this paper, we show that $\IAut(A)$ is also metabelian-by-(locally finite). In particular, $\IAut(A)$ has a normal subgroup $\G$ such that $\IAut(A)/\G$ is locally finite and $\G'$ is an abelian periodic subgroup whose all subgroups are normal in $\G$. In the case when $A$ is periodic,
 $\IAut(A)$ results to be  abelian-by-(locally finite) indeed, while in the general case it is not  even (locally nilpotent)-by-(locally finite). Moreover,  we provide further details about the structure of $\IAut(A)$ in some other cases for $A$.}


%%%%%%%%%%%%%%%%%%%%%%%%%%%%%%%%%%%%%%%%%%%%%%%%%%%%%%%%%%%%%%%%%
%\newtheorem{T}{\textsc{$\bullet$ \ }}
%\newtheorem{Teo}{\textsc{Theorem}}
%\newtheorem{thm}{Teorema}[section]%.
%%%%%%%%%%%%%%%%%

%\begin{document}

%\titlerunning{Short form of title
%% if too long for running head
%
%\author{Ulderico Dardano -  Silvana Rinauro}

%\date{}

%\maketitle

%%%%%%%%%%%%%%%%%%%%%%%%%%%%%%%%%%%%%%%%%%%%%%%%%%%%%%%%%%%%%%%%%%%%%%%%%%%%%%%%%55
%%%%%%%%%%%%%%%%%%%%%%%%%%%%%%%%%%%%%%%%%%%%%%%%%

\section{Introduction}

An endomorphism $\g$ of an abelian group $(A,+)$ is called  \emph{inertial endomorphism} if and only if  $( H+H\g)/H$ is finite for each subgroup $H$ (see \cite{DR2}, \cite{DGBSV}). An inertial endomorphism which is bijective is called {\it inertial automorphism}.
This definition can be seen as a generalization of the notion of  \emph{finitary automorphism} of $A$, that is, an automorphism $\g$ acting as the identity map on a  subgroup  of finite index in $A$ (see \cite{BS1}, \cite{W}). Since  $A$ is abelian, this condition is clearly equivalent to $A(\g-1)$ being finite. Note that we regard abelian groups as right modules over their endomorphism ring and  reserve the letter $A$ for abelian groups, which are {\it additively written}. 

The concept of {inertial endomorphism} of an abelian group $A$  may be used as a tool in the study of  \emph{inert} subgroups of possibly non-abelian groups (see \cite{DR5}, \cite{DGBSV2}). 
 Recall that a subgroup is called inert if it is
 {commensurable} with its conjugates (see \cite{BKS}, \cite{R}),
where  subgroups $H$ and $K$ are called {\it commensurable} if and only if $H\cap K$ has finite index in both $H$ and $K$.

In this paper we study the group
$\IAut(A)$ {\it generated} by all inertial automorphisms  of an abelian group $A$. Recall that  in \cite{DR1}  (resp. \cite{DR2}) we gave  
a description of inertial automorphisms (resp. endomorphisms) of an abelian group, while the ring of inertial endomorphisms of $A$ was  featured in \cite{DR4}. In  particular,  from \cite{DR2}, we have:\\  - { $\IAut(A)$ consists of products $\g_1\g_2^{-1}$ where $\g_1$ and $\g_2$ are both inertial automorphisms},\\ - { $\IAut(A)$ is locally (central-by-finite)},\\  
- { $\IAut(A)$ is
abelian modulo its subgroup $\FAut(A)$ of finitary automorphims}.\\ Recall that $\FAut(A)$ is known to be locally finite (\cite{W}).

Note that the above definitions of $\IAut(G)$ and $\FAut(G)$ make sense even if the underlying group $G$ is not abelian, and that $\FAut(G)\le \IAut(G)$ in any case. Also, in \cite{BS1} it has been shown that the group $\FAut(G)$ of finitary automorphisms of any group $G$ is both abelian-by-(locally finite) and (locally finite)-by-abelian.

\medskip
\noindent \textbf{Question}. \emph{Is $\IAut(A)$, the group  generated by all inertial automorphisms of an abelian group $A$, abelian-by-(locally finite)?}
\medskip

Our main results, in Sect. 3, can be summarized as follows:\\
- Theorem A  {gives a complete description of the group $\IAut(A)$, when $A$ is periodic};\\
- Corollary A {asserts that  the answer to our question is in the positive, when $A$ is periodic};\\ 
- Theorem B {asserts that, for any abelian group $A$,  the group  $\IAut(A)$ has a metabelian subgroup $\G$ such that $\IAut(A)/\G$ is locally finite and each subgroup of $\G'$ is normal in $\G$};\\ 
- Corollary B {asserts that the answer is in the negative for $A=\Z(p^\infty)\oplus\Q_p$};\\ 
- Theorem C {describes the group $\IAut(A)$ in some cases when $A$ is non-periodic}.

\medskip
In our investigations, we shall look for  abelian normal subgroups $\Sigma$ of $\IAut(A)$ such that  the automorphisms induced by 
$\IAut(A)$ via conjugation on $\Sigma$ are inertial, see Theorem C. Thus, in Sect. 4, 
 we highlight the role played by stability groups with respect to finitary and inertial automorphisms of $A$. In Sect. 5 we treat the case when $A$ is periodic  (Theorem A). In Sect. 6 we treat the remaining cases (Theorems B and C).

\bigskip 

{\small
For undefined terminology, notation and basic facts we refer to \cite{F} or \cite{Rm}. In particular, $\pi(n)$ denotes the set of prime divisors of $n\in\Z$. If $\pi$ is a set of primes, then  $A_\pi$, $T(A)$ and $D(A)$ denote the unique maximum $\pi$-subgroup, the torsion subgroup and the divisible subgroup of the abelian group $A$, respectively. By {exponent} $m=exp(A)$ of a $p$-group $A$ we mean 
the smallest $m$ such that $p^{m}A=0$, or $m=\infty$ if $A$ is unbounded.
 Furthermore, 
$\Q_p$ is the additive group of rational numbers whose denominator is a power of the prime $p$ and $\Z(p^\infty):=\Q_p/\Z$.
Also, $r_0(A)$ denotes the torsion-free rank of $A$, i.e.  the cardinality of a maximal $\Z$-independent subset of $A$.  

If $A_i\le A$ and $\g\in \Aut(A)$, as usual we   denote by $\g_{|A_i}$ the restriction of $\g$ to $A_i$. 
If  we have $A=A_1\plus A_2$ and $\g_{|A_i}=\g_i\in \Aut(A_i)$ (for $i=1,2$), we  write $\g=\g_1\oplus\g_2$.  

Commutators are calculated in the holomorph group $A\rtimes \Aut(A)$. Moreover, if $\p$ is an endomorphism of the additive abelian group $A$ and $a\in A$, we use the notation $[a,\p]:=a\p-a=a(\p-1)$.

%%%%%%%%%%%%%%%%%%%%%%%%%%%%%%%%%%%%

\section{Preliminaries}

 It is convenient to recall that an automorphism leaving every subgroup invariant is  usually called a {\it power automorphism}. Then 
 the group $\PAut(A)$ of  power automorphisms of an abelian group $A$ can be described as follows  (see \cite{Rm}). 

If $A$ is a $p$-group and $\alpha=\sum\limits_{i=0}^{\infty }\alpha_{i}p^{i}$ (with $0\le \alpha_{i}<p$) is an invertible  $p$-adic, we define, with an abuse of notation,  the power automorphism $\alpha$ (that we will also call \emph{multiplication} by $\alpha$) by setting $a\alpha:=(\sum\limits_{i=0}^{k-1 }\alpha_{i}p^{i}) a$, for any $a\in A$ of order $p^{k}$. In this way, we have defined an action  on $A$ of the group $\U_{p}$ of units of the ring of  $p$-adic integers, whose image is $\PAut(A)$. If $A$ has infinite exponent, then this action is faithful and $\PAut(A)$ is isomorphic to $\U_{p}$. Otherwise, if $e:=exp(A)<\infty$,  then the kernel of this action is $\{ \alpha\in \U_{p}\  |\ \alpha\equiv 1$ {\rm mod}\ $p^{e}\}$ and $\PAut(A)$ is isomorphic to the group of units of $\Z(p^{e})$. 

If $A$ is any periodic abelian group, then $\PAut(A)$ is the cartesian product of all the $\PAut(A_{p})$ where $A_{p}$ is the $p$-component of $A$.  If $A$ is non-periodic, then  $\PAut(A)=\{\pm1\}$.

 According to \cite{DR2}, an automorphism $\g$ is called an (invertible) {\it  multiplication}  of $A$ if and only if it  is a  power automorphism of $A$, if $A$ is periodic, or -when $A$ is non-periodic- there exist  coprime integers $m, n$ such that $(na)\g=ma$,  for each $a\in A$.  In the latter case, we have $mnA=A$ and $A_{\pi(mn)}=0$ and -with an abuse of notation-  we will write $\g=m/n$. We warn that  we are using the word ``multiplication" in a way different from \cite{F}.  
 Invertible multiplications of $A$ form a subgroup which is a central subgroup of $\Aut(A)$.

 If $r_0(A)<\infty$,  from \cite{DR2} we have that $\IAut(A)$ contains the group of all invertible {multiplications}. In this case,
{  $\IAut(A)$ consists of inertial automorphisms only}. Furthermore,  {$\IAut(A)$  is the kernel of the setwise action of $\Aut(A)$ on the quotient  of the lattice of the subgroups of $A$ with respect to commensurability} (which is a lattice congruence indeed, since $A$ is abelian).

{\ If  $r_0(A)=\infty$,  then the above kernel is the subgroup of $\IAut(A)$ consisting of so-called \emph{almost-power} automorphisms of $A$}, that is, automorphisms $\g$ such that every subgroup {contains} a $\g$-invariant subgroup of finite index. This group was introduced in \cite{FGN} to study generalized soluble groups in which subnormal subgroups are normal-by-finite (or core-finite, according to the terminology of \cite{BLNSW} and \cite{CKLRS}). 

\bigskip

We recall now some other facts that will be used in the sequel. They follow from Theorem 3 of \cite{DR1},  Proposition 2.2 and Theorem A of \cite{DR2}.

%%%%%%%%%%%%%%%%%%%%%%%%%%%%%%%%%%%%%%%%%%%%%%%%%%%%%%%

\bigskip \noindent 
\begin{lemma}\label{Recalls} 
Let $\g$ be an automorphism of an abelian group $A$.\\ 
\noindent 
{1)} If $r_{0}(A)=\infty$, then $\g$ is inertial if and only if there are a subgroup  $A_{0}$ of finite index in $A$ and an   integer $m$  such that $\g=m$ on   $A_{0}$.\\ 
\noindent 
{2)} If $0<r_{0}(A)<\infty$, then $\g$ is inertial if and only if there are  a torsion-free $\g$-finitely generated $\g$-subgroup $V$ such that $A/V$ is periodic,  a rational number $m/n$  (with $m$,$n$ coprime integers) such that $\g=m/n$ on $V$  and  $A_{\pi}$ is bounded, where $\pi:=\pi(mn)$. In particular $A/A_\pi$ is $\pi$-divisible.\\ 
\noindent  
{3)} If $A$ is periodic, then $\g$ inertial if and only if $\g$ is inertial on each $p$-component of $A$ and acts as  a power automorphism on all but finitely many of them.\\
{4)} If $A$ is a $p$-group, then $\g$ is inertial if and only if either $\g$  acts as an invertible multiplication  (that is as a power automorphism)  on a subgroup $A_{0}$ of finite index in $A$ or (critical case)  $0\ne D:=D(A)$ has finite rank, $A/D$ is infinite   bounded and there is a subgroup $A_{1}$ of finite index in $A$ such that $\g$ acts as invertible multiplication (by  possibly different $p$-adics) on both $A_1/D$ and $D$.\qed 
\end{lemma}

\medskip
\noindent For further instances
of inertial automorphisms see  Lemma \ref{semimoltiplicazioni}.

%%%%%%%%%%%%%%%%%%%%%%%%%%%%%%%%%%%%%%%%%%%%

\section{Main results}%%%%%%%%%%

Our first result is about periodic abelian groups.\bigskip

\noindent \textbf{Theorem A.} {\it
Let $A$ be a periodic abelian group. Then there is a subgroup $\Delta$ of $ \IAut(A)$ 
which is direct product of finite abelian groups and such that 
$$\IAut(A)=\PAut(A) \cdot \FAut(A) \cdot  \Delta$$\vskip-1mm
\noindent where $\Delta$ is trivial, 
if $A$ is reduced. 

 Moreover, there are a set $\pi$ of primes and  subgroups $\Sigma$, $\Psi$ of $\IAut(A)$ such that 
 $\Sigma$ is 
 an abelian $\pi'$-group  with bounded primary components
and  
$$\FAut(A) \cdot  \Delta=\FAut(A_\pi) \times   (\Sigma \rtimes\ {\Psi})  $$
 where the automorphims induced by $\Psi$ via conjugation on $\Sigma$ are inertial and  this  action is faithful. 
}

\smallskip
\noindent \textbf{Corollary A.} {\it If $A$ is a periodic abelian group, then $\IAut(A)$ is central-by-(locally finite)}.\medskip

\noindent With an abuse of notation, in Theorem A we regard $\FAut(A_\pi)$ as naturally embedded in $\FAut(A)$. For details in the case $A$ is a $p$-group see Proposition \ref{IAut(p-group)} below. 

\bigskip

In the next theorem we answer our question in the non-periodic case.
We reduce to study the subgroup  $\IAut_1(A)$ consisting of  {\it  inertial automorphisms of $A$ that act as the identity map on $A/T(A)$}. Actually, applying results from \cite{DR1}, we have that the above introduced group of  \emph{almost-power automorphisms} of $A$ is 
%\\ \centerline{$\IAut_1(A)\times\{\pm1\}$.} 
$$\IAut_1(A)\times\{\pm1\}.$$
In the case when $r_0(A)=\infty$, we have $\IAut_1(A)=\FAut(A)$ by Lemma \ref{Recalls}.(1). 

We will also consider a group  \emph{$Q(A)$} of particular inertial automorphisms of $A$, which is cointained in the center of $\Aut(A)$ and 
is naturally isomorphic to the multiplicative group of rational numbers generated by $-1$ and primes $p$ such that  $A/A_p$ is $p$-divisible and 
 $A_p$ is either bounded or finite according as $r_0(A)$ is finite or not  (see Lemma \ref{semimoltiplicazioni} below for details). 

\bigskip

\noindent {\bf Theorem B.} \emph{ Let $A$ be a non-periodic abelian group. Then there is a
subgroup $Q(A)$, which is  isomorphic to a multiplicative group of  rational numbers, such that
$$\IAut(A)=\IAut_1(A)\times Q(A)$$}
%%%%%%%%%%%%%%%%%%%%%%%%%%%%%
\emph{Moreover there is a normal subgroup $\G$ of $\IAut_1(A)$ such that:\\
i)\ \ $\IAut_1(A)/\G$ is locally finite;\\
ii)\ the derived subgroup $\G'$ of $\G$ is a periodic abelian group and each  subgroup of $\G'$ is normal in $\G$}

\medskip \noindent
\textbf{Corollary B.} \emph{ If $A$ is an abelian group, then $\IAut(A)$ is metabelian-by-(locally finite). However,  $\IAut(\Z(p^\infty)\oplus \Z)$ is not nilpotent-by-(locally finite)}.

\bigskip 
\noindent   Note that if $A$ is torsion-free, then $\IAut(A)=Q(A)$ is abelian, as in Theorem 2 of \cite{DR1}. Further, in the statement of Theorem B  one may take $\G$ to be the subgroup of $\IAut_1(A)$ consisting of inertial automorphisms acting by multiplication on $T(A)$. Unfortunately this subgroup need not  be nilpotent, as in Corollary B.  On the other hand, 
 groups with property $(ii)$ in Theorem B above have been studied under the name of KI-groups in a series of papers (see \cite{S-ZD}). 

The next theorem considers cases in which $A$ splits on its torsion subgroup. 
For details see Propositions \ref{IAut(AsplitsTbounded)} and \ref{IAut()}.

\medskip
\noindent \textbf{Theorem C.}\quad \emph{If $A$ is an abelian group with  $r_0(A)<\infty$ and either $T:=T(A)$ is bounded or $A/T$ is finitely generated, then
 there are subgroup $\Sigma$ and $\Gamma_1$ of $\IAut_1(A)$ such that  
$$\IAut_1(A)=\Sigma\rtimes \G_1$$
where  $\Sigma$ is  a periodic abelian group, $\Gamma_1\iso IAut(T)$ and the  automorphisms induced by  $\G_1$ via conjugation on $\Sigma$ are  inertial}

\medskip
When $A/T$ is not finitely generated, it may happen that $A$ has very few inertial automorphisms, since from Proposition \ref{few} and  Lemma \ref{Recalls}.(2) we have 
$$\IAut(\Z(p^\infty)\oplus \Q_{(p)} )=\{\pm1\}.$$
However, in the general case, the group  $\IAut_{1}(A)$ may be large, see 
 Remark \ref{controesempioAbelxLF}.

%%%%%%%%%%%%%%%%%%%%%%%%%%%%%%%%

%%%%%%%%%%%%%%%%%%%%%%%%%%%%%%%%%%%%%%%%%%%%%%%%
%
%
%
%%%%%%%%%%%%%%%%%%%%%%%%%%%%%%%%%%%%%%%%%%%%%%%%%%%%%%%

%%%%%%%%%%%%%%%%%%%%%%%%%%%%%%%%%%%%%%%%%%%%%%%%%%%%%%%%%%%%%%%%%%%%%%%

\section{Finitary automorphisms and Stability groups}

We state now some basic facts that  perhaps are already known (see also  \cite{CP}). 
If $X\le A$, we denote by $\St(A,X)$ the stability group of the series $A\ge X\ge 0$, that is, the set of $\g\in \Aut(A)$ such that $X\ge[A,\g]:=A(\g-1)$ and $[X,\g]=0$. When $X$ is a characteristic subgroup of  $A$, each $\g\in \Aut(A)$ acts via conjugation on the abelian normal subgroup $\Sigma:=\St(A,X)$ of $\Aut(A)$, according to the rule $\sigma\mapsto \g^{-1}\sigma\g=:\sigma^\g$ for each $\sigma\in\Sigma$. Similarly, $\g$ acts on
the additive group $\Hom(A/X,X)$ of homomorphisms $A/X\to X$ by a corresponding formula, i.e. $\p\mapsto \g_{|A/X}^{-1}\p\g$ where $\p\in \Hom(A/X,X)$ and $\g_{|A/X}$
denotes the group isomorphism  induced by $\g$ on $A/X$. With an abuse of notation,  we denote by $\sigma -1$ the well-defined homomorphism $\bar a\in A/X \mapsto a\sigma-a\in X$. \bigskip

\begin{Fact} \label{Fact} \it The map\ \ \ ${\cal H}:\sigma\in St(A,X) \mapsto(\sigma-1)\in \Hom(A/X,X)$\ \ \ 
is an isomorphism of (right) $\Aut(A)$-modules, that is, for each $\g\in \Aut(A)$ we have 
$$\sigma^{\g}=\g^{-1}(\sigma-1)\g+1$$\vskip-8mm\qed
\end{Fact}

By this argument we have two technical lemmas. For the first one see \cite{DF2} .

\begin{lemma}\label{LemmaDF2} Let $A$ be an abelian group, $\sigma,\g\in \Aut(A)$ and $m_{1}, m_{2}\in\Z$. If $\sigma$ stabilizes a series $0\le A_1\le A$,  where   $\g=m_{1}$ on $A_{1}$ and $\g^{-1}=m_{2}$ on $A/A_{1}$, then $\sigma^\g=\sigma^{m_1m_2}$.\qed
\end{lemma}

Our next lemma deals with the case when $A$ splits on $X$ and will be used several times.  In such a condition, once fixed a direct decomposition $A=X\oplus K$, we have an embedding $\Aut(K)\to \Aut(A)$ given by $\g\mapsto 1\oplus\g $.  Note that, if $\G\n \Aut(A)$, then one can consider $\St_{\G}(A,X):=\St(A,X)\cap \G$ which is $\Aut(A)$-isomorphic to a submodule of $\Hom(A/X,X)$. The proof of the lemma is straightforward.

\begin{lemma}\label{LemmaSemidir}
Let $A=X\oplus K$, where $X$ is a $\G$-subgroup, $\G\le \Aut(A)$,  $\zeta:A/X\leftrightarrow K$ the natural isomorphism, $\Sigma:=\St_\G(A,X)$,
$\G_1:=\{\g_{|X}\oplus1 |\ \g\in \G\}$ and
$\G_2:=\{1\oplus\zeta^{-1}\g_{|A/X}\zeta\ |\ \g\in \G\}$. Then:\\
1) if $\G_1\le \G$, then $\G=C_\G(X)\rtimes \G_1$ and $C_\G(A/X)=\Sigma\rtimes \G_1$;\\
2) if $\G_2\le \G$, then $\G=C_\G(A/X)\rtimes \G_2 $ and $C_\G(X)=\Sigma\rtimes \G_2$;\\
3) if $\sigma\in\Sigma$,\ $\g_1\in C_\G(A/X)$ and $\g_2\in C_\G(X)$, then
$$\sigma^{\g_1\g_2}=\g_2^{-1}(\sigma-1)\g_1+1.$$
In particular, if $\G_1\G_2\le\G$, then $\G=\Sigma\rtimes (\G_1\times\G_2)$.\qed
\end{lemma}

\begin{proposition}\label{LemmaFAut}  Let $A$ be an abelian group and  $T:=T(A)$. \\ 
1) If $r_0(A)<\infty$, then the automorphisms induced by $\FAut(A)$  via conjugation on $\St(A,T)$ are finitary;\\
2) If  $r_0(A)=\infty$ and the quotient $A/T$ is free abelian, then 
there is $\g\in\FAut(A)$ which induces via conjugation on $\St(A,T)$ a non-finitary automorphism, provided $\FAut(T)\ne1$.
\end{proposition}

\pf \ \ $1)$\quad Denote $\bar A=A/T$ and fix $\g\in \FAut(A)$. By  Fact \ref{Fact}, for each $\sigma\in \St(A,T)$ we have   $[\sigma,\g]{\cal H}=(\sigma^{-1}\sigma^\g){\cal H}=-(\sigma-1)+(\sigma^\g-1)=
-(\sigma-1)+
(\sigma-1)\g=(\sigma-1)(\g-1)=:\p_\sigma$. 
Thus we have to check that the set $\{ \p_\sigma\ | \  \sigma\in  \St(A,T)\}$ is  finite. For each 
$ \sigma$, we have that $im(\p_\sigma)\le im(\g-1)$ has finite order, say $n$.
On the other hand, $\ker(\p_\sigma)\ge n\bar A$ and $\bar A/n\bar A$ is finite since $\bar A$ has finite rank.

$2)$\quad  If $A=T\oplus K$, where $K$ is free abelian on the infinite $\Z$-basis $\{a_i\}$, take $\g_{0}\in \FAut(T)\setminus\{1\}$. Let $t\in T$ such that $t\g_{0}\ne t$ and $\g:=\g_0\oplus 1$. For each $i$ define $\sigma_i\in \St(A,T)$ by the rule $a_i(\sigma_i-1):=t$ and $a_j(\sigma_i-1):=0$ if $j\ne i$. Then there are infinitely many $[\sigma_i,\g]$, as $a_i\not\in \ker([\sigma_i,\g]{\cal H}) \ni a_j$ for each $i\ne j$.
\QED  

%%%%%%%%%%%%%%%%%%%%%%%%%%%%%%%%%%%%%%
Clearly, it may well happen that $\St(A,T)\not\le \FAut(A)$, as  in the case $A=\Z(p^\infty)\oplus\Q_p$. On the other hand, we do have $\St(A,X)\le \FAut(A)$, provided that one of the following holds:\\
- $A/X$ is bounded and $X$ has finite rank, as in Propositions \ref{FAut(critico)} and \ref{IAut(p-group)}.(2);  \\
-  $A/X$ has finite rank and $X$ is bounded, as in Proposition \ref{IAut(AsplitsTbounded)}; \\ 
- $A/X$ is finitely generated and $X$ is periodic, as in Proposition \ref{IAut()}. \\
Here we consider an instance of the first case with $X=D(A)$ and prove a proposition concerning finitary automorphisms.

%%%%%%%%%%%5
\begin{proposition}\label{FAut(critico)} Let $A$ be an abelian $p$-group such that $D:=D(A)$ has finite rank and $A/D$ is bounded. 
Then  $\Sigma:=\St(A,D)$ is a bounded abelian $p$-group and there is a subgroup $\Phi\iso \FAut(A/D)$ such that 
$$\FAut(A)=\Sigma\rtimes \Phi$$ 
 where the automorphisms induced by $\Phi$ via conjugation on $\Sigma$ are finitary and this action is faithful.\end{proposition}
%%%%%%

\pf First note that if $\sigma\in\Sigma$, then $[A,\sigma]=A(\sigma-1)$ is finite, since it is both finite rank and bounded. Hence $\sigma\in \FAut(A)$. Consider a decomposition $A=D\oplus B$ and  apply Lemma \ref{LemmaSemidir}, with $X=D$ and $\G=\FAut(A)=C_\G(X)$. Put $\Phi:=\G_2$. Then $\FAut(A)=\Sigma\rtimes \Phi$, as claimed.

Let $\g \in \Phi$. We have to show that  set $\{[\sigma,\g]\ |\ \sigma\in \Sigma\}$ is finite. As in Proposition  \ref{LemmaFAut}, we have   $[\sigma,\g]{\cal H}=(\sigma^{-1}\sigma^\g){\cal H}=(1-\sigma)+(\sigma^\g-1)=
-(\sigma-1)+\g^{-1}(\sigma-1)=(\g^{-1}-1)(\sigma-1)=:\p_\sigma$. Thus we have to count how many homomorphisms $\p_\sigma$ there are. On the one hand, $\ker(\p_\sigma)$ contains $\ker(\g^{-1}-1)$ which has finite index in $A/D$. On the other hand, the image of each $\p_\sigma$ is contained in the finite subgroup $D[p^m]$, where $p^m$ is a bound for $A/D$. Therefore, there are only finitely many $\p_\sigma$, once  $\g$ is fixed.

Let us check that the action is faithful. Let $1\ne\g \in \Phi$ and let $b\in B$ with maximal order and $b\ne b\g$. Then $B=\la b\ra \oplus B_0$ and we can write $b\g=nb+b_0$ with $n\in\Z, b_0\in B_0$.
If $b\ne nb$, then
 there is $\sigma\in\Sigma$ such that $B_0(\sigma-1)=0$ and $b(\sigma-1)=d$
where $d\in D$ has the same order as $b$. Thus, by Fact \ref{Fact},  $b\g(\sigma^\gamma-1)=b\g(\g^{-1}(\sigma-1))=d$, while $b\g( \sigma-1)=nd$. Therefore $\sigma^\gamma\ne\sigma$. Similarly, if
$b= nb$, then there is $\sigma\in\Sigma$ such that $b(\sigma-1)=0$ and $b_0(\sigma-1)=d_1$ of order $p$.
Then $b\g(\sigma^\gamma-1)=0$, while $b\g( \sigma-1)=d_{1}$ and again $\sigma^\gamma\ne\sigma$.
\QED

\begin{remark}\label{RemarkFC} In Proposition \ref{FAut(critico)}, $\Sigma$ need not  be  contained in the FC-center of $\FAut(A)$.
\end{remark}

\pf Write $A=D\oplus B_0$  where $D\iso\Z(p^\infty)$ and $B_0=\bigoplus_{i}\la b_{i}\ra\le B$ is infinite and homogeneous.
 Fix $\sigma\in\Sigma$ such that $b_1(\sigma-1)=d$, where $d$ is an element  of $D$ of order $p$, and $\sigma-1=0$ on $D\oplus\bigoplus_{j\ne 1}\la b_{j}\ra$.
 For each $i$ consider $\g_i\in \FAut(A)$ switching $b_{i}\leftrightarrow b_{1}$ and acting trivially on $D\oplus (\bigoplus_{j\not\in\{1, i\}}\la b_{j}\ra)$. Then $\sigma_{i}^\g=\g^{-1}_{i}(\sigma-1)+1$. Hence
 $b_i\sigma^{\g_{i}}=d+b_i$ and $b_j\sigma^{\g_{i}}=b_j$ for each $j\ne i$.
\QED

%%%%%%%%%%%%%%%%%%%%%%%%%%%%%%%%%%%%%%%%%%

Now an instance of a similar argument with $X=T(A)$

\begin{proposition}\label{IAut(Asplits)} 
 Let $A$ be an abelian group with  $r_0(A)<\infty$ such that $A/T$ is  finitely generated (resp. $T:=T(A)$ is bounded).  Then $\Sigma:=\St(A,T)$ is a periodic (resp. bounded) abelian group and 
 there is a subgroup  $\Phi_1\iso \FAut(T)$ such that 
$$\FAut(A)=\Sigma\rtimes \Phi_1$$  
 where  
$\Phi_1$ induces via conjugation on $\Sigma$  finitary automorphims.

If $A/T\ne 0$ is finitely generated, then this action is faithful, while if $A=\Z_{12} \oplus \Q_{(2)}$ it is not.
\end{proposition}

\pf In any case, we can write $A=T\plus K$ where $r:=r_0(K)<\infty$. Recall that $\Sigma\iso \Hom(A/T,T)$.
Note that $\Sigma\le \FAut(A)$. In fact, if $\sigma\in\Sigma$, then $\sigma-1\in \Hom(A/T,T)$ and $A(\sigma-1)$ is an abelian grougp which is both finitely generated and periodic (resp. finite rank and bounded). Hence $A(\sigma-1)$ is finite that is  $\sigma\in\FAut(A)$.

Clearly  $\Phi_1:=\{\p\oplus1\ |\ \p\in \FAut(T)\}\iso\FAut(T)$ and $\Phi_1\le\FAut(A)$.  By Lemma  \ref{LemmaSemidir}.(1) we have that  $\FAut(A)=\Sigma\rtimes \Phi_1$.  By Proposition   \ref{LemmaFAut}, $\Phi_1$  induces via conjugation on $\Sigma$  finitary automorphisms.

If $A/T$ is  finitely generated, then $\Sigma\iso \Hom(A/T,T)$ is a periodic abelian group which is naturally isomorphic to the direct sum of $r$ copies of $T$ as a right $\Aut(A)$-module. Therefore the action of $\Phi_1$ on $\Sigma$ is faithful. Finally, if  $A=\Z_{12} \oplus \Q_{(2)}$, we have $\Phi_1\iso\U\Z_{12}$ and $\Sigma\iso\Z_3$; hence the action is not faithful.
\QED

%%%%%%%%%%%%%%%%%%%%%%%%%%%%%%%%%%%%%%%%%%%%%%%%%%%%%%%%%%%%%%%%%%%%%%%%%%

\section{The group $\IAut(A)$, when $A$ is periodic}

To give a detailed description of $\IAut(A)$ when $A$ a $p$-group, let us introduce some terminology.  By \emph{essential exponent} $e=eexp(A)$ of $A$ we mean
the smallest $e$ such that $p^{e}A$ is finite, or $e=\infty$ if $A$ is unbounded. In the former case, this is equivalent to saying that $A=A_0\oplus A_1\oplus A_2$ where  $A_0$ is finite, $exp(A_1)<e\le exp(A_0)$ and $A_2$ is the sum of infinitely many cyclic groups of order $p^e$.
In \cite{DR1} we called \emph{critical} a $p$-group of type $A=B\oplus D$ with $B$ infinite but bounded and $D\ne 0$ divisible with finite rank (see  Lemma \ref{Recalls}.(4)). Critical groups will be a tool to describe $\IAut(A)$ when $A$ is periodic.

\begin{proposition}\label{IAut(p-group)}  Let $A$ be an abelian $p$-group and $D:=D(A)$.

\smallskip\noindent 
 1) If $A$ is non-critical, then $\IAut(A)= \PAut(A) \cdot \FAut(A)$
 where $\PAut(A)\cap \FAut(A)$ is either trivial or cyclic of order $p^{m-e}$, according as $A$ is unbounded 
or $m:=exp(A)<\infty$ and $e:=eexp(A$).

\smallskip\noindent  
2) \ If $A=D\oplus B$ is critical, let  $\D:=\{1\oplus n \,|\, n\in\Z\setminus p\Z\}$, $\Phi:=\{1\oplus \p_{0} \,|\, \p_{0}\in\FAut(B)\}$ and $\Psi:=\{1\oplus \g_{0} \,|\, \g_{0}\in\IAut(B)\}$, then $$\IAut(A)= \PAut(A) \times (\FAut A\cdot \D).$$
Moreover\ \  $\FAut A\cdot \D=C_{\IAut(A)}(D)= \Sigma\rtimes\ \Psi$, where   $\FAut(A)=\Sigma\rtimes\Phi$ and\\
\phantom{xx}
 i)\ \ $\Sigma:=\St(A,D)$ is an infinite abelian $p$-group, $exp(\Sigma)=exp(B)=:m'<\infty$ and $eexp(\Sigma)=eexp(B)=:e'$;\\
\phantom{xx}
ii)\ \  $\Psi=\Phi\Delta\iso\IAut(B)$ where $[\Phi,\Delta]=1$ and $\Psi$ induces via conjugation  on $\Sigma$ inertial automorphisms and  this action is faithful;
\\
\phantom{xx} iii)\ \ $\D\iso\PAut(B)\iso \U(\Z(p^{m'}))$, each  $\delta_n:=1\oplus n\in\D$ acts via conjugation  on $\Sigma$ as the multiplication by $n$ and 
$\FAut(A)\cap\D$ has order $p^{m'-e'}$;\\
\phantom{xx}  iv)\ \ $\Phi\iso\FAut(B)$ and $\Phi$
 induces via conjugation on $\Sigma$  finitary automorphisms.
\end{proposition} 

%\smallskip
\pf 
Let $\g\in \G:=\IAut(A)$. 

$1)$ If $A$ is non-critical, then, according to  Lemma \ref{Recalls}.(4), there exist a $p$-adic $\alpha$ and a  subgroup $A_{0}$ of finite index in $A$ such that $\g_{|A_{0}}=\alpha$. Thus $\g^{-1}\alpha$ acts on $A_{0}$ as the identity map, that is, $\g^{-1}\alpha\in \FAut(A)$. Hence $\IAut(A)= \PAut(A) \cdot \FAut(A)$.
Further,  if the $p$-adic number $\beta$ is in $\PAut(A)\cap \FAut(A)$, then $\beta$ is trivial on a
 subgroup $B$ of finite index in $A$. Therefore $\beta=1$, provided
$exp(A)=\infty$. Otherwise,
$exp(B)\ge e$ and $\beta\equiv 1\ mod\ p^e$. Thus there are at most $p^{m-e}$ choices for such a $\beta$. On the other hand, each $p$-adic number $\beta\equiv 1\ mod\ p^e$ is  finitary.

\medskip

$2)$ Let $A=D\oplus B$ be critical. By   Lemma \ref{Recalls}.(4) there exists an invertible $p$-adic $\alpha$ such that $\g_{|D}=\alpha$. Thus $\g_1:=\g\alpha^{-1}\in C_\G(D)$. Clearly, $\PAut(A)\cap C_\G(D)=1$, so that
$\IAut(A)= \PAut(A) \times C_\G(D)$.  

 Again by  Lemma \ref{Recalls}.(4), $\g_1$ acts by multiplication by an integer $n$ on a  subgroup of finite index in $A[p^{m'}]$ where $A[p^{m'}]\ge B$. Therefore, if $\delta_n:=1\oplus n\in \D$ with respect to $A=D\oplus B$, we have $\g_1\delta_n^{-1}\in \FAut(A)$.  Hence $C_\G(D)=\FAut(A)\cdot\D$.

It is routine to verify that  $(i)$ holds, since $\Sigma:=\St(A,D)\iso Hom(B,D)$. By Proposition \ref{FAut(critico)}, $(iv)$ holds as well. By Lemma \ref{LemmaSemidir} (with $X:=D$, $K:=B$  and so $\G_2=\Psi$), we have $C_\G(D)= \Sigma\rtimes\ \Psi$ as stated in $(2)$. Then, applying part $(1)$ of the statement to $B$, we have $\Psi=\Delta\Phi$ and $[\Phi,\Delta]=1$ as in $(ii)$. Moreover, $\FAut(A)\cap\D$ has order $p^{m'-e'}$. 

By Lemma \ref{LemmaDF2}, we have that $\Delta$ acts on $\Sigma$ as in $(iii)$. Thus the whole $\Psi=\Phi\Delta$ acts via conjugation on $\Sigma$ inducing inertial automorphisms and $(ii)$ holds.

It remains to show that $\Psi$ acts faithfully on $\Sigma$. Let $\p\delta_{n}\in C_{\Psi}(\Sigma)$
with $\p\in\Phi$ and $\delta_{n}:=1\oplus n\in\Delta$. On the one hand, $\delta_{n}$ acts via conjugation  on $\Sigma$ as the multiplication by $n$ by $(iii)$. On the other hand,  $\delta_n$  is finitary on $\Sigma$ by $(iv)$. Since $eexp(\Sigma)=eexp(B)$ by $(i)$, then multiplication by $n$ is finitary on $B$. Thus $\delta_{n}\in C_{\Phi}(\Sigma)=1$ by Proposition \ref{FAut(critico)}.  \phantom{xxxx}
\QED

We have seen that, if $A$ is a $p$-group, then $\IAut(A)$ is central-by-(locally finite). If $A$ is a critical $p$-group, one can ask whether there is an abelian normal subgroup $\Lambda$ of $\IAut(A)$ such that $\IAut(A)=\Lambda\cdot \FAut(A)$. The answer is in the negative, as in the following remark. Fisrt we state an easy  lemma.

\begin{lemma}\label{decomposition} If $B_0$ is a subgroup of finite index in a bounded abelian group $B$, then there are subgroups $B_1$ and $B_2$  such that $B_2$ is finite, $B_1\ge B_0$ and 
$B=B_1\oplus B_2.$
\end{lemma}

\pf Clearly there is a finite subgroup $F$ such that $B=B_0+F$. Since $B_0$ is separable and $B_0\cap F$ is finite, then there is a finite subgroup $B_3\ge B_0\cap F$ such that $B_0=B_1\oplus B_3$ for some $B_1\le B_0$. Fix $B_1$
 and $B_2:=B_3+F$. On the one hand $B_1+B_2=B_1+B_3+F=B_0+F=B$. On the other hand, by Dedekind law, $B_1\cap B_2=
 B_1\cap (B_3+F)= B_1\cap (B_0\cap (B_3+F))=B_1\cap (B_3+(B_0\cap F))=B_1\cap B_3 =0$.\QED

\begin{remark}\label{remark} If $A$ is a critical $p$-group (with $p\ne2$) and $\Lambda\n \IAut(A)$ is such that $C_\G(D)=\FAut(A)\cdot \Lambda$, then $\Lambda$ is neither finite nor locally nilpotent.
\end{remark}

\pf We use the same notation as in Proposition \ref{IAut(p-group)}. Let $n\in \N$ be  a primitive root of $1$  mod $p^{m'}$ and consider  $\delta:=1\oplus n\in\D$ with respect to $A=D\oplus B$. Since 
$\D\le C_\G(D)=\FAut(A)\cdot \Lambda$,  then we can write  with $\p\in \FAut(A)$ and $\lambda\in\Lambda$. Hence $\delta=\lambda=n$ on some subgroup $B_0$ of  finite index in $B$. By Lemma \ref{decomposition}, $B=B_1\oplus B_2$ with $B_1\le B_0$ and $B_2$ finite. Put $A_1:=D+B_1$ and note that ${\lambda}_{|A_{1}}=1\oplus n$ with respect to $A_1=D\oplus B_1$. 

It is sufficient to show that $\la{\lambda}\ra^{\G_{1}}$ is infinite and not locally nilpotent, where $\G_{1}$ is the group of (inertial) automorphisms of $A$ of type $ \g_1\oplus1$ with respect to  $A=A_{1}  \oplus B_{2}$, with $\g_1\in \IAut (A_{1})$. Thus we may assume $A_{1}=A$ and $\G:=\G_1$. Then multiplication by $n$ is in $\Lambda$ and $\Lambda=\D^\G$.

{We claim that $\D^\G=\Sigma \rtimes \D$}. In fact, by  Proposition \ref{IAut(p-group)} we have that $\D\iso \U(\Z_{p^{m'}})$ acts faithfully by multiplication on the infinite abelian $p$-group $ \Sigma$  of exponent ${m'}$   and then $\Sigma=[\Sigma,\D]$ and $\D^\G=\Sigma\D$, as claimed. Thus  $\D^\G$ is not locally nilpotent, since the action of $\D$ on $\Sigma$ is fixed-point-free.\QED

\noindent {\bf Proof of Theorem A}. By  Lemma \ref{Recalls}.(3),  $\IAut(A)$ may be identified with $\PAut(A)\cdot Dr_p\,\, \IAut(A_p)$.
Apply Proposition \ref{IAut(p-group)} to each $A_p$. Let $\pi$ be the set of primes $p$ for which $A_p$ is not critical. If $p\in\pi$, we have 
$\IAut(A_p)=\PAut(A_{p})\cdot \FAut (A_p)$. Otherwise, for each $p\not\in \pi$,   
there are subgroups $\D_p,\Sigma_p,{\Psi_p}$ corresponding to $\D,\Sigma,{\Psi}$ in Proposition \ref{IAut(p-group)} 
 such that $\IAut(A_p)=\PAut(A_{p})\cdot \FAut (A_p)\cdot \D_p$ and  $\FAut (A_p)\cdot\D_p = \Sigma_p\rtimes\, {\Psi_p}$. Now  
it is routine to verify that  the  statement follows by setting 
$\D:=Dr_{p\not\in\pi}\ \D_p$, $\Sigma:=
 Dr_{p\not\in\pi}\ \Sigma_p$, $\Psi:=Dr_{p\not\in\pi}\ \Psi_p$,
 and recalling that $Dr_{p}\ \FAut(A_p)=\FAut(Dr_{p}\ A_p)$.  \qed
\bigskip

Remark that, in Theorem A,  when we consider the action of the above $\Psi$  on the $p$-component $\Sigma_p$ of $\Sigma$ we are concerned with subgroups of $\IAut(\Sigma_p)=\PAut(\Sigma_p)\cdot \FAut(\Sigma_p)$, where $\Sigma_p$ is a bounded abelian $p$-group and $\PAut(\Sigma_p)$ is finite abelian.

%%%%%%%%%%%%%%%%%%%%%%%%%%%%%%%%%%%%%%%%%%%%%%%%%%%%%%%%%%%
%%%%%%%%%%%%%%%%%%%%%%%%%%%%%%%%%%%%%%%%%%%%%%%%%%%%%%%%%%%%%%%%%%%%%%

%%%%%%%%%%%%%%%%%%%%%%%%%%%%%%%%%%%%%%%%%%%%%%%%%%%%%%%%%%%

\section{The group $\IAut(A)$, when $A$ is non-periodic}

%%%%%%%%%%%%%%%%%%%%%%%%%%%%%%%%%%%%%%%%%%%%%%%%%5

To prove Theorem B  we point out the existence of some inertial automorphisms of a particular type.

\begin{lemma}\label{semimoltiplicazioni} 
 Let $A$ be a non-periodic abelian group
and $\pi_*(A)$ be the set of primes such that  $A/A_p$ is $p$-divisible and one of the following holds:\\
- $A_p$ is finite,\\ 
- $r_0(A)$ is finite and $A_p$ is  bounded.

Then, for each $p\in\pi_*(A)$, there is a unique $C^{(p)}$ such that $A=A_p\oplus C^{(p)}$ and the automorphism $\g_{(p)}:=1\oplus p$ (with respect to this decomposition) is inertial.

Moreover, the subgroup $Q(A):=\la \g_{(p)}\ |\ p\in\pi_*(A)\ra\times \{\pm1\}$ is a central subgroup of $\IAut(A)$ which is isomorphic to a multiplicative group of rational numbers.  Furthermore, $\IAut_1(A)\cap Q(A)=1$. 
\end{lemma}

\pf The proof of the first part of the statement, concerning $C^{(p)}$, is routine.  Let us show that $\g_{(p)}$ is inertial. For each $H\le A$ we have $H+H\g_{(p)}\le H+A_p$.  If $A_p$ is finite, then $(H+H\g_{(p)})/H$ is finite as well and $\g_{(p)}$ is inertial. If $r_0(A)$ is finite  and $A_p$ is  bounded, let  
 $V$ be the $\la\g_{(p)}\ra$-closure of a free subgroup of $C$ with maximal rank. Then 
$V$ is torsion-free, $C/V$ is a $p'$-group and $\g_{(p)}$ acts as a power automorphism  on $A/V$. Thus we apply Lemma \ref{Recalls}.(2) and  deduce that $\g_{(p)}$ is inertial. The remaining part of the statement follows straightforward.\QED

\noindent {\bf Proof of Theorem B}. Let $\g_{(p)}$  and $Q:=Q(A)$ as  in Lemma \ref{semimoltiplicazioni}.  

We first consider the case when $r_0(A)=\infty$.
Let $\g\in \IAut(A)$. Then, by Corollary B in \cite{DR2}, we have   $\g=\g_1\g_2^{-1}$ with $\g_1$,$\g_2$ inertial. Further, by  Lemma \ref{Recalls}.(1), there is a subgroup $A_0$ with finite index in $A$ such that we have $\g_{|A_0}=m/n=p_1^{s_1}\cdots p_t^{s_t}\in\Q$ ($m,n$ coprime, $p_{i}$ prime, $s_{i}\in \Z$). Also $\IAut_1(A)=\FAut(A)$ and $\g=m/n$ on $A/T$ as well. If $m/n=1$, then $\g\in \FAut(A)$. If $m/n=-1$, put $\g_0:=-1\in Q$. Otherwise,
since $\g$ is invertible, $mA_0=A_0=nA_0$.
Then for each $p_i\in\pi:=\pi(mn)$,
 the $p_i$-component of $A$ is finite and $A/T$ is $p_i$-divisible.  Consider then   $\g_0:= \g_{(p_1)}^{s_1}\cdots \g_{(p_t)}^{s_t}\in Q$. In both cases, $\g\g_0^{-1}=1$ on $A_0/(A_0)_\pi$ hence
 $\g\g_0^{-1}\in \FAut(A)$. Thus $\IAut(A)=\IAut_1(A)\times Q(A)$. 
Moreover, $(i)$ and $(ii)$ are true with $\G=1$, since $\IAut_1(A)=\FAut(A)$ is locally finite.

Let then $r_0(A)<\infty$ and $\g\in \IAut(A)$. By Corollary B in \cite{DR2}  $\g$ is inertial. By  Lemma \ref{Recalls}.(2), we have that 
$\g=m/n=p_1^{s_1}\cdots p_t^{s_t}\in\Q$ ($m,n$ coprime, $p_{i}$ prime, $s_{i}\in \Z$) on $A/T$. We also have that, for each $p_i\in\pi:=\pi(mn)$, the group  $A/T$ is $p_i$-divisible and $A_{p_i}$ is bounded. Consider 
  $\g_0:= \g_{(p_1)}^{s_1}\cdots \g_{(p_t)}^{s_t}\in Q$, Clearly  $\g_0=m/n$ on $A/T$. Thus $\g\g_0^{-1}$ acts trivially on $A/T$ and $\IAut(A)=\IAut_1(A)\times Q(A)$, as stated.

Let $\G$  be  the preimage of $\PAut(T)$ under the canonical homomorphism $\IAut_1(A)\mapsto \IAut(T)$. Now $(i)$ holds, since
$\IAut_1(A)/\G$ is locally finite by Theorem A.
To check $(ii)$ consider that the derived subgroup $\G'$ of $\G$ stabilizes the series $0\le T\le A$ and therefore is abelian. Moreover, by Theorem B in \cite{DR2}, the subgroup  $\G'$ consists of finitary automorphisms. Thus  $\G'$  is torsion and  $(ii)$  holds by Lemma \ref{LemmaDF2}. \qed

%%%%%%%%%%%%%%%%%%%%%%%%%%%%%%%%%%%
\medskip

Let us see that there are groups $A$ with few inertial automorphisms even if $r_{0}(A)<\infty$ and that the canonical homomorphism $\IAut_{1}(A)\to \IAut(T)$ need not be surjective.

%%%%%%%%%%%%%%%%%%%%%%%%%%%%%%%%%

\begin{proposition}\label{few} Let $A$ be a $\pi$-divisible non-periodic abelian group, where $\pi$ is a set of primes. If $T:=T(A)$ is a $\pi$-group, then $\IAut_1(A)=1$. 
\end{proposition}
%%%%%%%%%%%%%%%%%%%%%%%%%%%%%%%%%%%%%%%%%%%%

%\vskip2mm

\pf  If $r_0(A)=\infty$ then $\IAut_1(A)=\FAut(A)$. Moreover, if $\g\in \FAut(A)$ then $A(\g-1)$ is a  finite $\pi$-group. Then $A/\ker(\g-1)$ is such. Hence $A=\ker(\g-1)$ and $\FAut(A)=1$.

If $r_0(A)<\infty$, by  Lemma \ref{Recalls}.(2) we have $\g=1$ on some free abelian subgroup $V\le A$ such that $A/V$ is periodic. Moreover, the $\pi$-component $B/V$ of $A/V$ is  divisible.   Then,  by  Lemma \ref{Recalls}, part (3) and (4), we have  that 
$\g$ is a multiplication on $B/V$. Furthermore, the group $B/(V+T)$ is $\pi$-divisible and has non-trivial $p$-component for each $p\in\pi$, since $(V+T)/T\iso V$ is free abelian. Thus from $\g=1$ on $B/T$ it follows that $\g=1$ on 
 $\g=1$. Hence $\g$ stabilizes the series $0\le V\le B$. However $\Hom(B/V,V)=0$. Then $\g=1$ on $B$. Therefore $\g-1$ induces a  homomorphism $A/B\rightarrow T$ which is necessarily $0$ since $A/B$ is a $\pi'$-group. Thus $\g=1$ on the whole group $A$. \QED

%%%%%%%%%%%%%%%%%%%%%%%%%%%%%%%%%%%%%%%%%%%%%%%%%%%%%%%%%%%%%%%%%%%%%%%%%%%%%
\vskip2mm

\noindent {\bf Proof of Theorem C}. It follows from the next propositions which considers cases in which $\IAut_1(A)$ splits on $\Sigma:=\St(A,T)$.\QED

%%%%%%%%%%%%%%%%%%%%%%%%%%%%%%%%%%%%%%%%%%%%%%

\begin{proposition}\label{IAut(AsplitsTbounded)} 
 Let $A$ be an abelian group  and  $T:=T(A)$.

If   $r_0(A)<\infty$ and $T$ is bounded, then $\Sigma:=\St(A,T)$ is a  bounded abelian group and there is a subgroup $\G_1$ of $\IAut_{1}(A)$ such that  $\G_1\iso \IAut(T)$ and 
 $$\IAut_{1}(A)=\Sigma\rtimes\G_1$$
where  $\G_1$ induces via conjugation on $\Sigma$ inertial automorphisms.
\end{proposition}

\pf We can write $A=T\plus K$ where $r:=r_0(K)<\infty$. Note that the group $\Sigma\iso \Hom(A/T,T)$ is a periodic abelian group which  is bounded as $T$.

Clearly  $\G_1:=\{\g \oplus1\ |\ \g\in \IAut(T)\}\iso\IAut(T)$. If $\g\in \IAut(T)$, then $\g\oplus 1$ (with respect to $T\plus K$) is inertial by  Lemma \ref{Recalls}.(4), and so $\G_1\le\IAut_{1}(A)$. Thus we may  apply  Lemma \ref{LemmaSemidir} with $\G:= \IAut_1(A)$. We obtain
 $\IAut_{1}(A)=\Sigma\rtimes\G_1$, as claimed.

  By Proposition \ref{IAut(p-group)}, we have
$\IAut(T)=\FAut(T)\cdot  \PAut(T)$.
Hence $\G_1=\Phi_{1}\Delta_1$ where $\Phi_1:=\{\p\oplus1\ |\ \p\in \FAut(T)\}\iso\FAut(T)$ acts conjugation on $\Sigma$ by means of finitary automorphisms, by Proposition \ref{IAut(Asplits)} and $\Delta_1:=\{\delta\oplus1\ |\ \delta\in \PAut(T)\}\iso \PAut(T)$  acts via conjugation on $\Sigma$ by means of multiplications,  by Lemma \ref{LemmaDF2}.
Therefore the whole $\G_1$ induces via conjugation  on $\Sigma$  inertial automorphisms.
\QED

\noindent We notice that the action of $\G_1$  on $\Sigma$ in Proposition \ref{IAut(AsplitsTbounded)}  need not be faithful, as already seen in Proposition \ref{IAut(Asplits)}.

\begin{proposition}\label{IAut()} 
 Let $A$ be a non periodic abelian group  and  $T:=T(A)$.

If $A/T$ is  finitely generated, then $\Sigma:=\St(A,T)$ is a periodic abelian  group
 and 
there is a subgroup $\G_1$ of $\IAut_1(A)$ such that  $\G_1\iso \IAut(T)$ and 
 $$\IAut_{1}(A)=\Sigma\rtimes\G_1$$
where  $\G_1$ induces via conjugation on $\Sigma$ inertial automorphisms and this action is faithful. 

If in addition $T$ is unbounded, then $\IAut_{1}(A)$ is not nilpotent-by-(locally finite). Further, if $A_{2'}$ is unbounded, then $\IAut_{1}(A)$ is not even (locally nilpotent)-by-(locally finite).
\end{proposition}

\pf As in the proof of Propositon \ref{IAut(AsplitsTbounded)}, we can write $A=T\plus K$ where $K$ is finitely generated. The group $\Sigma\iso \Hom(A/T,T)$ is a periodic abelian group which is isomorphic to the direct sum $\oplus_r T$ of $r:=r_0(A)>0$ copies of $T$ as a right $\Aut(A)$-module. 

Clearly  $\G_1:=\{\g \oplus1\ |\ \g\in \IAut(T)\}\iso\IAut(T)$. If $\g\in \IAut(T)$, then  $\g\oplus 1$ (with respect to $T\plus K$) is inertial by  Lemma \ref{Recalls}.(4). Hence  $\G_1\le\IAut_{1}(A)$. Thus we may  apply  Lemma \ref{LemmaSemidir} with $\G:= \IAut_1(A)$,  and we obtain
 $\IAut_{1}(A)=\Sigma\rtimes\G_1.$ 
 
Let us investigate now the action of $\G_1$ via conjugation on $\Sigma$. Assume first that $T$ is a $p$-group. 
 Let $\g\in \IAut(T)$. By  Proposition \ref{IAut(p-group)}, $\g=\g_{0}\p$, where  $\p\in \FAut(T)$ and either $\g_{0}\in \PAut(T)$  or $T$ is a critical $p$-group and $\g_{0}$ induces multiplications on both $D(T)$ and $T/D(T)$. Recall that $\Sigma$ is $\Aut(A)$-isomorphic to $\oplus_r T$.
In the former case,
 that is if $\g_{0}\in \PAut(T)$, then $\g_0\oplus 1$ acts via conjugation on $\Sigma$ as a power automorphism (that is a multiplication).
In the latter case, $\Sigma$ is critical as well and $\g_{0}\oplus 1$ induces invertible multiplications on both $D(\Sigma)$ and
$\Sigma/D(\Sigma)$. Thus $ \g_{0}\oplus 1$ acts via conjugation on  $\Sigma$ as an inertial automorphism of $\Sigma$,  by  Lemma \ref{Recalls}.(4).  In both cases, by Proposition  \ref{IAut(Asplits)}, $\p$ acts via conjugation on $\Sigma$ as a finitary automorphism. 
Hence $\g\oplus 1$ acts via conjugation on  $\Sigma$ as an inertial automorphism. 

In the general case, when $T$ is any periodic group and $\g\in \IAut(T)$, then $\g\oplus 1$ (with respect to $T\plus K$)  acts via conjugation as an inertial automorphism on all primary components $\Sigma_p$ of $\Sigma$, by what we have seen above and the fact that $\Sigma_p\iso \Hom(A/T,A_p)$. Similarly, since $\g\oplus 1$ acts as a multiplication on all but finitely many primary components $A_p$ of $A$, it acts the same way on all but finitely many $\Sigma_p$. Thus $\g\oplus 1$ is inertial on $\Sigma$ by  Lemma \ref{Recalls}.(3).

It is clear that the action via conjugation of $\G_1$ on $\Sigma$ is faithful
 as the standard action of $\G_1$  on $T$ is such.

To prove the last part of the statement, note that in the case when $T$ is unbounded, then there exists a non-periodic multiplication $\alpha$ of $T$. Note that the automorphism $\mu:=\alpha\oplus 1$ (with respect to $T\plus K$) belongs to $\G_1$. If, by the way of contradiction, $\la\Sigma,\mu\ra$ is nilpotent-by-(locally finite), then there
is $s\in\Z\setminus\{0\}$ such that $\la\Sigma,\mu^s\ra$ is nilpotent, so there is $n\in\N$ such that 
$[\Sigma,_n\mu^s]=0$, and hence $0=\Sigma^{(\mu^s-1)^n}= \Sigma^{(\alpha^s-1)^n}$. This is  a contradiction, since $\Sigma$ is unbounded as $T$ is.

Finally, if $A_{2'}$ is unbounded, then $\Sigma_{2'}$ is unbounded as well. Let $\alpha$ be a non-periodic multiplication of $A_{2'}$. Then, $\mu:=\alpha\oplus 1\oplus 1$   with respect to  $A=A_{2'}\oplus A_2\plus K$ acts as non-periodic multiplication (by $\alpha$) of $\Sigma_{2'}$ acting fixed-point-free on a primary component. Thus $\mu$ (and any non-trivial power of $\mu$ as well) does not belong to the locally nilpotent radical $R$ of $\IAut_1(A)$. Therefore $\IAut_{1}(A)/R$ is not locally finite.
\QED

%%%

\vskip-1mm
%%%%%%%%%%%%%%%%%%%%%%%%%%%%%%%%%%%%%%%%%%%%%%%%%%%%55
Finally, we note that, despite the above propositions, in the general case the group  $\IAut_{1}(A)$ may be large.
\vskip-1mm

\begin{remark}\label{controesempioAbelxLF}
There exists an abelian group $A$ with $r_0(A)=1$ and $A_p \iso\Z(p)$ for each prime $p$ such that $\IAut(A)=\IAut_1(A)\times\{\pm1\}$,   $\IAut_1(A)=\Sigma\cdot \FAut(A)$, where $\Sigma:=\St_{\IAut(A)}(A,T(A))\not\le \FAut(A)$, $\Sigma\iso \prod_p\Z(p)$ and $\IAut_1(A)/\FAut(A)\iso \Sigma/T(\Sigma)$ is a  divisible torsion-free abelian group with cardinality $2^{\aleph_0}$.

Moreover  any  element of $\IAut_1(A)$ induces on $T$ a finitary automorphism.
\end{remark}
\vskip-1mm
\pf   As in Proposition A in \cite{DR2}, we consider the group  $G:=B\oplus C$ where
$B:=\prod_p\langle b_{p}\rangle$,
$C:=\prod_p\langle c_{p}\rangle$, and $b_p$, $c_p$ have order $p$, $p^2$ resp.  and $p$ ranges over all primes. Consider then the (aperiodic) element $v:=(b_p+pc_p)_p\in G$ and $V:=\la v\ra$. We have that for each prime $p$ there is an element $d_{(p)}\in G$ such that $pd_{(p)}=v-b_{p}$. Let 
$A:=V+\langle d_{(p)}|\ p\ \rangle$.
 Then $A/T\iso \la 1/p\ |\ p\ \ra\le \Q$, since $A/T$ has torsion free rank $1$ and $v+T$ has $p$-height $1$ for each $p$. Thus $T=T(B)\iso\bigoplus_p \Z(p)$ and 
 the $p$-component of $A/V$ is generated by $d_{(p)}+V$ and has order $p^2$, since  $pd_{(p)}=v-b_{p}$.

Then $\Sigma\iso \prod_p\Z(p)$ and $\Sigma\cap F\Aut(A)=T(\Sigma)$, hence $\Sigma\not\le F\Aut(A)$. Moreover $A=\la d_{(p)}\ra+V$, where $V=\la v\ra$ is infinite cyclic and $A_p=\la b_p\ra$ has order $p$. Also $\Aut(A/T)=\{\pm1\}$ and $\IAut(A)=\IAut_1(A)\times\{\pm1\}$.

We claim that \emph{if  $\g\in \IAut_1(A)$ induces  on $T$ a finitary automorphism, then $\g\in \Sigma\cdot \FAut(A)$}. In fact, $T\g$ is finite, so it is a $\pi$-component of $A$ for some finite $\pi$. Thus $\g\g_{0}^{-1}\in \Sigma$, where $\g_{0}:=\g_{|A_{\pi}}\oplus 1$ with respect to $A=A_{\pi}\oplus K$ and clearly $\g_{0}\in \FAut(A)$.

 Finally we prove the last part of the statement, from which it follows  $\IAut_1(A)=\Sigma\cdot \FAut(A)$. Let  $\g\in \IAut_1(A)$ and $\p:=\g-1$.  Since $A\p\le T$,  there exists an integer $n\not=0$ such that $(nv)\p=0$. 
We prove that   $T\p\subseteq A_{\pi(n)}$, which is finite.  For any prime $p$, on the one hand, $nd_{(p)}$ is a $p$-element modulo $\la nv\ra\le \ker \p$, hence
$(nd_{(p)})\p\in A_p$, that implies $(pnd_{(p)})\p=p(nd_{(p)})\p=0$. On the other hand,  $(pnd_{(p)})\p=n(v-b_p)\p=-n(b_{p})\p$. Hence, if $p\not\in (n)\p$, then $A_{p}\p=0$.
\QED

%%%%%%%%%%%%%%%%%%%%%%%%%%%%%%%%%%%%%%%%%%%%%%%%%
%%%%%%%%%%%%%%%%%%%%%%%%%%%%%%%%%%%%%%%%%%%%%%%%%%%%%%%%%%%%%%%%%%%%%%%

%%%%%%%%%%%%%%%%%%%%%%%%%%%%%%%%%%%%%%%%%%%%%%%%%%%%%%%%%%%%%%%%%%%%%%%%%%

%ultima_sezione_dopo_\end{document}

%%%%%%%%%%%%%%%%%%%%%%%%%%%%%%%%%%%%%%%%%%%%%%%%%%%%%%%%%%%%%%%%%%%%%%%%%%%%%%%%%%%%%%%%%%%%%%%

%%%%%%%%%%%%%%%%%%%%%%%%%%%%%%%%%%%%%%%%%%%%%%%%%%%%%%%%%%%%%%%%%%%%%%%%%%%%%%%%%%%%
%\newpage
\small \vskip-1cm
{

}

%%%%%%%%%%%%%%%%%%%%%%%%%%%%%%%%%%%%%
%%%%%%%%%%%%%%%%%%%%%%%%%%%%%%%%%%

\bigskip
\noindent  Ulderico Dardano,  Dipartimento di Matematica e
Applicazioni ``R.Caccioppoli'', Universit\`a di Napoli ``Federico
II'', Via Cintia - Monte S. Angelo, I-80126 Napoli, Italy.
\texttt{dardano@unina.it}

\medskip
\noindent  Silvana Rinauro,
Dipartimento di Matematica, Informatica ed Economia, Universit\`a della
Basilicata, Via dell'Ateneo Lucano 10,
I-85100 Potenza, Italy. \texttt{silvana.rinauro@unibas.it}

%%%%%%%%%%%%%%%%%%%%%%%%%%%%%%%%%%%
%%%%%%%%%%%%%%%%%%%%%%%%%%%%%%%%%%%%
%%%%%%%

\end{document}